\documentclass[12pt]{article}
\usepackage{amssymb}

\newtheorem{theo}{Theorem}

\newtheorem{prop}{Proposition}

\newtheorem{lem}{Lemma}

\newtheorem{rem}{Remark}

\newtheorem{defi}{Definition}

\newtheorem{exemp}{Example}
\begin{document}
\author{ Aziz Ikemakhen \thanks{
Facult\'{e} des sciences et techniques, B.P. 549, Gueliz,
Marrakech, Maroc. \hspace{4cm} \ \  E-mail:
ikemakhen@fstg-marrakech.ac.ma } }
\title{ Parallel Spinors on Pseudo-Riemannian $Spin^c$ Manifolds}
\date{  }
\maketitle

\begin{abstract}
 We describe, by their holonomy groups,  all complete  simply connected  irreducible
 non-locally symmetric pseudo-Riemannian $Spin^c$
  manifolds which admit parallel spinors. So
  we generalize the Riemannian $Spin^c$ case (\cite{Mor}) and the pseudo-Riemannian
  $Spin$ one (\cite{BK}).
\end{abstract}
\section{Introduction}
In (\cite{Mor}),  Moroainu described all complete simply connected
 Riemannian $Spin^c$ manifolds admitting parallel spinors.
Precisely, he showed the following result:
\begin{theo}
 A complete simply
connected $Spin^c$ Riemannian  manifolds $(M,g)$ admits a parallel
spinor if and only if it is isometric to the Riemannian product
$(M_1,g_1)\times (M_2,g_2) $ of a complete simply connected
k$\ddot{a}$hler manifold $(M_1,g_1)$ and a complete simply
connected $Spin$ manifold $(M_2,g_2) $ admitting a parallel
spinor.  The $Spin^c$ structure of $(M,g)$ is then the product of
the canonical $Spin^c$ structure of $(M_1,g_1)$ and the $Spin$
structure of $(M_2,g_2)$.
\end{theo}
In (\cite{BK}), Baum and Kath  characterized, by their holonomy
group, all simply connected  irreducible non-locally symmetric
pseudo-Riemannian $Spin$ manifolds admitting parallel spinors.
Precisely, they proved the following result:
\begin{theo}
Let $(M,g)$ be a simply connected  irreducible
 non-locally symmetric pseudo-Riemannian $Spin$ manifold of
 dimension $n=p+q$ and  signature $(p,q)$. We denote by $N$  the dimension of the
 space of parallel spinors on $M$.
 Then   $(M,g)$ admits
  a  parallel spinors if and only if the
 holonomy group $H$ of $M$ is (up to conjugacy in $O(p,q)$)
 one  in the following list :

$$
\begin{tabular}{||c|c||}
\hline \hline
  Holonomy group & N \\
  \hline
  \hline
  $SU(p',q') \subset SO(2p',2q')$ & 2 \\
  \hline
  $Sp(p',q') \subset SO(4p',4q')$ & $ p' + q' + 1$ \\
  \hline
  $ G_2 \subset SO(7)$ & 1 \\
  \hline
  $ G'_{2(2)} \subset SO(4,3)$ & 1 \\
  \hline
  $ G^{\mathbb{C}}_2 \subset SO(7)$ & 2 \\
  \hline
  $ Spin(7)\subset SO(8)$ & 1 \\
  \hline
  $ Spin(4,3) \subset SO(8, 8)$ & 1 \\
  \hline
  $Spin(7,\mathbb{C})\subset SO(8, 8)$ &  1 \\
  \hline
\end{tabular}
$$
$$(table \;\;1)$$
\end{theo}
Our aim is to generalize this result for the simply connected
irreducible non-locally symmetric pseudo-Riemannian $Spin^c$
manifolds. More  precisely, we show that:
\begin{theo} Let $(M,g)$ be a simply connected  irreducible
 non-locally symmetric pseudo-Riemannian  manifold of
 dimension $n=p+q$ and  signature $(p,q)$. Then
 the following conditions are equivalent\\
 (i) $(M,g)$ is a  $Spin^c$ manifold which admits
  a  parallel spinor,  \\
 (ii) Either $(M,g)$ is a  $Spin$ manifold which admit a parallel
 spinor,
 or $(M,g)$ is a k$\ddot{a}$hler not special k$\ddot{a}$hler
 manifold,\\
 (iii) The holonomy group $H$ of $(M,g)$ is (up to conjugacy in
$O(p,q)$)
 one  in table 1 or $H = U(p',q')$, $p = 2p'$ and $q = 2q'$.\\
 For $H = U(p',q')$ the dimension of the space of parallel spinors on $M$ is 2.
\end{theo}
This theorem is a contribution to the resolution of the following
problem :\\
\textbf{(P)} What are the possible holonomy groups of simply
connected
 pseudo-Riemannian $Spin^c$ manifolds which admit parallel
 spinors?\\
 Some partial answers to this problem have been given by M. Wang  for
 the Riemannian $Spin$ case (\cite{Wa}), by H. Baum and I. Kath  for the
 irreducible pseudo-Riemannian $Spin$ one (\cite{BK}) , by  Th. Leistner
for the  Lorentzian $Spin$ one (\cite{Le}, \cite{L2}), by  A.
Moroainu for the Riemannian $Spin^c$ one (Theorem 1), and by
author for the totally reducible pseudo-Riemannian $Spin$ one and
for Lorentzian $Spin^c$ one (\cite{I1}, \cite{I2} ). The problem
remains open even though big progress
have been made. \\
The proof of Theorem 3 is based on Theorem 2 and the technique
used by Moroianu to prove Theorem 1 that we adapted to
the pseudo-Riemannian case.\\
In paragraph 2 of this paper we define the group $Spin^{c}(p,q)$
and its spin representation. We also define the
$Spin^{c}$-structure on pseudo-Riemannian manifolds and its
associated spinor bundle. In paragraph 2 we give an algebraic
characterization to the pseudo-Riemannian $Spin^c$ manifolds which
admit parallel spinors and  we prove Theorem 3.
\section{ Spinor representations and  $Spin^{c}$- bundles }
\subsection{   $Spin^{c}(p,q)$ groups }
Let $<.,.>_{p,q}$ be the ordinary scalar product of signature
$(p,q)$ on $\mathbb{R}^m$ $(m=p+q)$. Let $Cl_{p,q}$ be the
Clifford algebra of $\mathbb{R}^{p,q}:=(\mathbb{R}^m
,<.,.>_{p,q})$ and $\mathbb{C}l_{p,q}$ its complexification. We
denote by $\cdot$ the Clifford multiplication of
$\mathbb{C}l_{p,q}$. $\mathbb{C}l_{p,q}$ contains the groups
$$\mathbb{S}^1:=\{z \in \mathbb{C};  \parallel z
\parallel=1  \}$$ and
$$Spin(p,q):=\{X_1 \cdot...\cdot X_{2k};  \;\;
<X_i,X_i>_{p,q} \; = \pm1;\;\; k \geq 0  \}.$$ Since $\mathbb{S}^1
\cap Spin(p,q)=\{-1,1\}$, we define the group $Spin^{c}(p,q)$ by
$$ Spin^{c}(p,q)= Spin(p,q)\cdot \mathbb{S}^1 \diagup
\{-1,1\} = Spin(p,q)\times_{\mathbb{Z}_2} \mathbb{S}^1 .$$
Consequently, the elements of $ Spin^{c}(p,q)$ are the
 classes $[g,z]$ of pairs $(g,z) \\ \in Spin(p,q)\times
\mathbb{S}^1 $, under  the equivalence relation $(g,z)\sim
(-g,-z)$. The following suites are exact (see \cite{LM}):
$$
\displaystyle 1\rightarrow \mathbb{Z}_2 \rightarrow Spin(p,q)
\stackrel {\lambda}{\longrightarrow} SO(p,q)\rightarrow 1
$$

$$
\displaystyle 1\rightarrow \mathbb{Z}_2 \rightarrow Spin^{c}(p,q)
\stackrel {\xi}{\longrightarrow} SO(p,q)\times \mathbb{S}^1
\rightarrow 1,
$$
where $\lambda(g)(x)= g \cdot x \cdot g^{-1}$ for $ x\in
\mathbb{R}^m $  and $ \xi ([g,z])=
(\lambda(g),z^2). $\\
Let $(e_i)_{1\leq i\leq m }$ be an  orthonormal basis of
$\mathbb{R}^{p,q}$ ( $<e_i,e_j> = \varepsilon_i \delta_{ij },$
$\varepsilon_i = -1 $ for $1\leq i\leq p $  and  $\varepsilon_i =
-1 $  for  $1+p\leq i\leq m $ ). The Lie algebras of $Spin(p,q)$
and $Spin^{c}(p,q)$ are respectively
$$ spin(p,q):= \{ e_i \cdot e_j \; ; 1 \leq i < j \leq m \} $$
and
$$
 spin^{c}(p,q):= spin(p,q) \oplus \texttt{i}\mathbb{R}.
 $$
The derivative  of  $\xi$ is a  Lie algebra isomorphism and it is
given by
$$
\xi_\ast ( e_i \cdot e_j , \texttt{i}t)= ( 2 E_{ij} , 2\texttt{i}t
),
$$
 where $ E_{ij}= -\varepsilon_j D_{ij}+ \varepsilon_i D_{ji}$ and $
 D_{ij}$ is the standard basis of  $gl(m,\mathbb{R})$ with the
 $(i,j)$-component equal  1 and all other zero.

\subsection{$Spin^{c}$ bundles}

In this paper, we will use the following isomorphisms:\\ Let $ U=
\left(
\begin{array}{cc}
 0 &  i   \\
 i &  0
\end{array}
\right), \;  V= \left(
\begin{array}{cc}
 0 &  -1   \\
 1 &  0
\end{array}
\right), \; E= \left(
\begin{array}{cc}
1 &  0   \\
 0 &  1
\end{array}
\right) \; , \; T= \left(
\begin{array}{cc}
 -1 &  0   \\
 0 &  1
\end{array}
\right) $, and $\mathbb{C}(2^n)$ the complex algebra consisting of
$2^n \times 2^n$-matrices.
\\
 In case $m =p + q = 2n$ is even, we define (see \cite{BK})
$ \Phi _{p,q} : \mathbb{C}l_{p,q} \rightarrow \mathbb{C}(2^n) $ by
$$ \Phi _{p,q}
(e_{2j-1})   =  \tau_{2j-1}  E \otimes ... \otimes E \otimes U
\otimes T \otimes ...\otimes T
 $$
\begin{equation}  \label{ eqn : 1}
 \;\;\;\; \;\; \Phi _{p,q} (e_{2j})   =   \tau_{2j}  E \otimes ... \otimes E \otimes V \otimes
  \underbrace{T \otimes ... \otimes T}_{\mbox{(j-1)-times}}.
 \end{equation}
Where $\tau_j =\mathrm{i} $ if $\varepsilon_j  = -1$ and $\tau_j
=1$ if
$\varepsilon_j  = 1.$\\
In case  $m =2n +1 $ is odd, $ \Phi _{p,q} :  \mathbb{C}l_{p,q}
\rightarrow \mathbb{C}(2^n) \oplus \mathbb{C}(2^n)$  is defined by

 $$
 \Phi _{p,q}(e_k)   = (\Phi _{p,q-1}(e_k),  \Phi
 _{p,q-1}(e_k)),\;\;
 k= 1, ..., m-1;
 $$
 \begin{equation}  \label{ eqn : 2}
\Phi_{p,q}(e_m)=(\texttt{i}T\otimes...\otimes T , -\texttt{i}
T\otimes...\otimes T ).
\end{equation}
This yields representations of the spin group $Spin (p,q)$ in case
$m$ even by restriction and in case $m$ odd by restriction and
projection onto the first component.
 the  module space of $Spin(p,q)$- representation is
  $ \Delta_{p, q}= \mathbb{C}^{2^n}$.
The   Clifford multiplication is defined by
 \begin{equation}  \label{ eqn : 3}
X \cdot u := \Phi_{p,q}(X) (u) \;\;\; \mbox{ for} \;\;\;  X \in
\mathbb{C}^m \;\;\; \mbox{and}  \;\;\;  u \in \Delta_{p, q}.
 \end{equation}
A usual basis  of $  \Delta_{p, q}$ is the following :
$u(\nu_n,...,\nu_1):=u(\nu_n) \otimes...\otimes u(\nu_1);
 \;\;   \nu_j =\pm1,$
 where \\

$u(1)=  \left(%
\begin{array}{c}
  1\\
  0
  \end{array}
  \right)   \;\; \hbox{and} \;\;
  u(-1)=  \left(%
\begin{array}{c}
  0\\
  1
  \end{array}
  \right ) \in \mathbb{C}^2.$ \\
  The spin representation of the group $Spin(p,q)$ extends to a
   $Spin^{c}(p,q)$- representation by :

\begin{equation}\label{4}
 \Phi_{p,q}([g,z])(v)= z \; \Phi_{p,q}(g) (v),
\end{equation}
for $v \in \Delta_{p,q}$ and $ [g,z] \in Spin^{c}(p,q)$. Therefore
 $\Delta_{p,q}$ becomes the  module space of
 $Spin^{c}(p,q)$- representation (see \cite{Frid}).\\
There exists a hermitian inner product $<.,.>_{\Delta}$ on the
spinor module $\Delta_{p,q}$ defined by
$$
<v,w>_{\Delta}:=  \texttt{i}^{\frac{p(p-1)}{2}}(e_1 \cdot ...\cdot
e_p\cdot v,w ) ; \;\; \hbox{for} \;\;  v, w  \in \Delta_{p,q},
$$
where $ (z,z')=\sum_{i=1}^{2^n} z_i \cdot \overline{z'_i}$ is the
standard  hermitian  product on  $\mathbb{C}^{2^n}$.\\
$<.,.>_{\Delta}$ satisfies the following properties :
\begin{equation}  \label{ eqn : 5}
<X\cdot v,w>_{\Delta} = (-1)^{p+1} < v,X\cdot w>_{\Delta},
\end{equation}
for $X \in \mathbb{C}^m$.
\subsection{Spinor bundles}

Let $(M,g)$\ be a connected pseudo- Riemannian manifold of signature $%
(p,q)$. And let  $P_{SO(p,q)}$\ denote the bundle of  oriented
positively frames on M.
\begin{defi}
A structure $Spin^{c}$ on  $(M,g)$ is the data of  a
$\mathbb{S}^1$-principal bundle  $P_L$ over $M$ and a
$\xi$-reduction $(P_{Spin^c(p,q)},\Lambda_{p,q})$ of the product
$(SO(p,q)\times \mathbb{S}^1)$-principal bundle  $P_{SO(p,q)}
\times P_L$.
i.e. $\Lambda: P_{Spin^c(p,q)} \rightarrow (P \times P_L)$ is a 2-fold covering verifying:\\
i) $P_{Spin^c(p,q)}$ is a $Spin^c(p,q)$-principal bundle over $M$,\\
ii) $\forall u \in Q$,  $\forall a \in Spin^\mathbb{c}(p,q)$,
$$ \Lambda(ua)= \Lambda(a)\xi(a). $$
What means, the following diagram commutes :
\[
\begin{array}{ccccc}
P_{Spin^c(p,q)}\times Spin^c(p,q) & \longrightarrow &
P_{Spin^c(p,q)} & \stackrel{\pi }{\searrow }
&  \\
\Lambda \otimes \xi \downarrow &  & \Lambda \downarrow &  & M. \\
(P_{SO(p,q)} \times P_L)\times (SO(p,q)\times \mathbb{S}^1) &
\longrightarrow & P_{SO(p,q)} \times P_L
 & \stackrel{\pi ^{\prime }}{%
\nearrow } &
\end{array}
\]
\end{defi}
Not that $(M,g)$ carries a $Spin^c$-structure if and only if the
second Stiefl-Whitney class of $M$, $w_2(M)$ is the mod reduction
of an integral class (\cite{LM}).
\begin{exemp} Every pseudo-Riemannian  $Spin$ manifold is
canonically a $Spin^c$ manifold. The $Spin^c$- manifold is
obtained as
$$
P_{Spin^c(p,q)} = P_{Spin(p,q)}\times_{\mathbb{Z}_2}\mathbb{S}^1,
$$
where $P_{Spin(p,q)}$ is the Spin-bundle and $\mathbb{Z}_2$ acts
diagonally by $(-1,-1)$.
\end{exemp}
\begin{exemp} Any  irreducible pseudo-Riemannian  k$\ddot{a}$hler manifold is
canonically a $Spin^c$ manifold.\\
Indeed The holonomy group $H$ of $(M,g)$ is $U(p',q')$, where
$(p,q)= (2p',2q')$ is the signature of $(M,g)$ . Then
 $P_{SO(p,q)}$ is reduced to $U(p',q')$-principal bundle
 $P_{U(p',q')}$. Moreover there exists an
 $<.,.>_{p,q}$-orthogonal almost complex structure  $J$
 which we can imbed $U(p',q')$ in
 $SO(p,q)$ by
 $$
\begin{array}{ccc}
  i : U(p',q') & \hookrightarrow & SO(p,q) \\
  A+ \texttt{i}B =((a_{kl})_{1\leq k ,l \leq m} + \texttt{i}(b_{kl})_{1\leq k ,l \leq m}) &
   \rightarrow & \left(%
   \begin{array}{c} \left(%
\begin{array}{cc}
  a_{kl} & b_{kl} \\
 -b_{kl} & a_{kl} \\
\end{array}%
\right)
\end{array}%
\right)_{1\leq k ,l \leq m}. \\
\end{array},
$$
and $( e_k  , J e_k )_{k=1,...,p'+q'}$ is an orthogonal basis of
$\mathbb{R}^{p,q}$. \\
 We consider the homomorphism
$$
\begin{array}{cccc}
  \alpha: & U(p',q') & \hookrightarrow & SO(p,q)\times \mathbb{S}^1 \\
   & C & \rightarrow & (i(C),\det(C)) \\
\end{array}
$$
Since the proper values of every element $C \in U(p',q')$ is in
$\mathbb{S}^1$ and
$$ \cos 2 \theta + \varepsilon_k \sin 2\theta \;\; e_k \cdot J e_k
= \varepsilon_k (\cos \theta \;\; e_k  +   \sin \theta \;\; J e_k
)\cdot (\sin \theta \;\; e_k  -   \cos \theta \; \; J e_k ),
$$

where $ \varepsilon_k = \; <e_k,e_k>_{p,q} $, then the following
homomorphism
$$
\begin{array}{cccc}
 \widetilde{\alpha}: & U(p',q') & \hookrightarrow & Spin^c(p,q)\times \mathbb{S}^1 \\
   & C & \rightarrow & \displaystyle \prod_{k=1}^m (
    \cos 2 \theta_k + \varepsilon_k \sin 2\theta_k \; \; e_k \cdot Je_k) \times
   \displaystyle e^{ \frac{\texttt{i}}{2} \sum \theta_k},\\

\end{array}
$$
is well defined, where $\displaystyle e^{i \theta_k}$ are the
proper values of $C$. And it is easy to verifies that

$$ \xi \circ \widetilde{\alpha} = \alpha.$$
Consequently,
$$
P_{Spin^c(p,q)} =
P_{U(p',q')}\times_{\widetilde{\alpha}}Spin^c(p,q).
$$
\end{exemp}
Now, let  denote by $S := P_{Spin(p,q)}
\times_{\Phi_{p,q}}\Delta_{p,q}$ the spinor bundle associated to
the $Spin^c$-structure $P_{Spin(p,q)}$, by $L :=P_L
\times_{i}\mathbb{C} = P_{Spin(p,q)}
\times_{\mathbb{S}^1}\mathbb{C}$ the complex line bundle
associated to the auxiliary bundle
  $P_L$,  where
 $$i : \mathbb{S}^1 \rightarrow
GL(\mathbb{C}).$$
The  Clifford multiplication given by (\ref{ eqn
: 3}) defines a Clifford multiplication on $S$:
$$
\begin{array}{ccc}
  TM \otimes S = (P_{Spin(p,q)} \times_{\Phi_{p,q}}\mathbb{R}^m )\otimes (P_{Spin(p,q)}
\times_{\Phi_{p,q}}\Delta^{\pm}_{p,q} )& \rightarrow & S \\
 ( X\otimes \psi ) = [q,x]\otimes [q,v]& \rightarrow & [q,x\cdot v]=:X\cdot\psi. \\
\end{array}
$$
Since the   scalar product $<.,.>_{\Delta}$ is $Spin_0
^c(p,q)$-invariant, it defines  a scalar product on $S$ by :
$$ <\psi,\psi_1>_{\Delta}= <v,v_1>_{\Delta},\;\; \hbox{for} \;\; \psi = [q,v]
\;\; \hbox{and} \;\; \psi_1 = [q,v_1] \in \Gamma(S).
$$
According to (\ref{ eqn : 5}), it is then easy to verify that
\begin{equation}  \label{ eqn : 6}
 <X\cdot \psi,\psi_1>_{\Delta} = (-1)^{p+1} <\psi,X\cdot
 \psi_1>_{\Delta},
\end{equation}
for $ X \in \Gamma(M)$ and   $\psi, \psi_1 \in \Gamma(S)$.\\
 Now,
as in the  Riemannian case ( see \cite{Frid}),
 if $(M, g)$ is a $Spin^{c}$ pseudo- Riemannian manifold, every connection form $A:
TP_L \rightarrow \mathrm{i} \mathbb{R}$  on the $\mathbb{S}^1$-
bundle $P_L$ defines ( together with the Levi-Civita $D$ of $(M,
g)$ ) a covariant derivative  $\nabla^A$ on the spinor bundle $S$,
called
the spinor derivative  associated to  $(M,g, S, L, A)$.\\
Henceforth, a $Spin^{c}$- pseudo- Riemannian manifold will be the
data of a set \\ $(M,g, S, P_L, A)$, where $(M,g)$ is an oriented
connected pseudo- Riemannian manifold, $S$ is a $Spin^{c}$
structure, $L$ is the complex line bundle associated to the
auxiliary bundle of $S$ and $A$ is a connection form on $P_L$.
Using (6) and by the same proof in the Riemannian case ( see
\cite{Frid}), we conclude that
\begin{prop} $ \forall \; X, Y \in \Gamma(M) $  and
$\forall \; \psi, \psi_1 \in \Gamma(S)$,
\begin{equation}  \label{ eqn : 7}
  \nabla^A _Y  (X\cdot
\psi) =X\cdot \nabla^A_Y  (\psi) + D_Y X \cdot \psi.
\end{equation}
\begin{equation}  \label{ eqn : 8}
 X< \psi,\psi_1>_{\Delta} =  <\nabla^A_X \psi,
 \psi_1>_{\Delta}+ <\psi,\nabla^A_X \psi_1>_{\Delta}.
\end{equation}
\end{prop}
Let us denoted by $F_A := \mathrm{i} w $ the  curvature form of
$A$, seen as an imaginary-valued  2-form on $M$,  by $R$ and $Ric$
respectively  the curvature
  and the Ricci tensor of $(M,g)$ and by $R^A$ the curvature tensor of $\nabla^A$.
  Like  Riemannian case ( see
  \cite{Frid}), if we put $
A(X):=X\lrcorner \omega $ we have
\begin{prop} For $q=(e_1,...,e_m)$ a local  section  of \ \ $P_{Spin(p,q)}$,\\
$ \forall \; X, Y \in \Gamma(M) $ \hbox{and} $\forall \; \psi \in
\Gamma(S)$,
\begin{equation}  \label{ eqn : 9}
R^A(X,Y)\psi = \frac{1}{2}\sum _{1\leq i<j\leq m}\varepsilon_i
\varepsilon_j \; g(R(X,Y) e_i , e_j ) e_i\cdot e_j\cdot \psi +
\mathrm{i} \frac{1}{2} \; \omega(X,Y)\cdot \psi.
\end{equation}
\begin{equation}  \label{ eqn : 10}
\sum _{1\leq i\leq m}\varepsilon_i \; R^A(X,e_i)\psi = -
\frac{1}{2} Ric(X)\cdot \psi + \mathrm{i} \frac{1}{2} A(X)
\cdot\psi.
\end{equation}
\end{prop}

\begin{rem}  According to Example 1, if $(M,g)$ is
 $Spin$
 then it is $Spin^{c}$. Moreover,
  the  auxiliary bundle $P_L$ is trivial and then there exists a
global  section  $ \sigma : M \rightarrow P_L$. We choose the
connection  defined by $A$ to be flat, and we denote $ \nabla^A $
by $\nabla .$  Conversely, if the auxiliary bundle $P_L$ of a
$Spin^{c}$-structure is trivial, it is canonically identified with
a
 $Spin$-structure. Moreover, if the  connection $A$ is flat,
by this identification, $\nabla^A$ corresponds to the
 covariant derivative on the spinor bundle.
\end{rem}

\section{parallel Spinors}
\subsection{algebraic characterization }

It is well know that there exists a bijection between  the space
$\mathcal{PS}$ of all parallel spinors on $(M,g)$ and the space
$$ V_{\widetilde{H}} = \{v \in
\Delta_{p,q}
 ; \;\; \Phi_{p,q}(\widetilde{H})(v):= \widetilde{H}\cdot v = v \}
 $$
of all fixed spinors of $\Delta_{p,q}$ with respect to the
holonomy group  $\widetilde{H}$ of the  connection $\nabla^A$. If
$(M,g)$ is supposed  simply connected, then $\mathcal{PS}$ is in
bijection with $$ V_{\widetilde{\mathcal{H}}}= \{v \in
\Delta_{p,q}
 ; \;\; \widetilde{\mathcal{H}}\cdot v = 0 \},$$
where $\widetilde{\mathcal{H}}$ is  the Lie algebra of
$\widetilde{H}$. With the notations  introduced in subsection 2.1,
for $B \in spin (p,q)$ and $t \in \mathbb{R} $, we have

 \begin{equation}  \label{ eqn : 11}
\xi^{-1}_{\ast} (B,  \mathrm{i} t) = ( \lambda^{-1}_{\ast} (B),
\frac{1}{2}\mathrm{i}t )= ( \frac{1}{4}\sum_{i=1}^m \varepsilon_i
\; e_i\cdot B(e_i) , \frac{1}{2}\mathrm{i}t ).
\end{equation}
Moreover, $\xi (\widetilde{H})= H \times H_A$, where $H$ is the
holonomy group of $(M,g)$ and $H_A$ the one of  $A$. $H_A =
\mathbb{S}^1$, if $A$ is flat and $H_A = \{1\}$ otherwise. Then

$$ V_{\widetilde{\mathcal{H}}}= \{v \in \Delta_{p,q}
 ; \;\; \xi_\ast^{-1}(\mathcal{H}\oplus \mathcal{H}_A)\cdot v=0
 \},
 $$
where  $\mathcal{H}$ is the   Lie algebra of  $H$ and
$\mathcal{H}_A $  the one of $H_A$.  $\phi_{p,q}$ is linear, then
if we differentiate the relation (4)  we get :
$$ \phi_{p,q}(B, \mathrm{i} t)(v) = \mathrm{i} t v + \phi_{p,q}(B)(v) =
\mathrm{i} t v + B \cdot v .$$  According to
 (11) we have ,
$$ \phi_{p,q}( \xi_\ast^{-1}(B, it))(v) =
 \frac{1}{2} \mathrm{i} t v + \lambda^{-1}_{\ast} (B)\cdot v.
$$
Therefore   $(M,g)$ admits a  parallel spinor  if and only if
there exists  $0 \neq v \in \Delta_{p,q}$ such that
\begin{equation}
\left\{%
\begin{array}{ll}
   \mathcal{H} \cdot v := \displaystyle \lambda_\ast^{-1}(\mathcal{H})\cdot v =
    \mathcal{H}_A \; v , \\
   \mathcal{H}_A = \{0\}  \;\; \hbox{or} \;\; \mathrm{i} \mathbb{R}. \\
\end{array}%
\right. \label{12}
\end{equation}

\subsection{  Proof  of  Theorem 3 }

Let $(M,g, S, P_L, A)$ be a $Spin^{c}$ structure where  $(M , g)$
is  a simply connected  irreducible
 non-locally symmetric pseudo-Riemannian  manifold of
 dimension $n=p+q$ and  signature $(p,q)$,
 which admits a non trivial parallel spinor
$\psi$. \\ We consider the two  distributions $T$ and $E$ defined
by
$$  T_x := \{ X \in T_xM;  \;\; X\cdot \psi = 0 \},
$$
$$
 E_x = \{ X \in T_xM; \; \exists \; Y \in T_xM; \;\; X\cdot \psi = \texttt{i}
 Y\cdot\psi \},
 $$
for $x \in M$.  Since $\psi$ is parallel, By (7), $T$ and $E$ are
parallel. Since  $T$ is isotropic and the manifold $(M,g)$  is
supposed irreducible, by the holonomy principe, we have
\begin{equation}\label{13}
   T=0.
\end{equation}
Now denote by $F$ the image of the Ricci tensor:
$$ F_x := \{ Ric(X);  \;\; X \in T_xM \}.$$
Since $\psi$ is parallel, (10) shows that
\begin{equation}\label{14}
   Ric(X)\cdot \psi = \texttt{i} A(X)\cdot \psi.
\end{equation}
Then $F \subset E$. Consequently, from (13), we have
$$  E^\perp  \subset F^\perp = \{ Y \in TM;\; Ric(Y)=0
 \} = \{ Y \in TM;\; A(Y )  = 0 \}.$$
 $(M,g)$  is supposed irreducible, by the holonomy
principe, $E= 0$ or $E= TM$. \\
If $E= 0$, then $F=0$. This gives $Ric=0$ and $A=0$. According to
Remark 1, $(M,g)$ is Spin and $\psi$ is a parallel spinor on
$M$.\\
If $E= TM$, we have a $(1,1)$-tensor $J$ definite by
\begin{equation}\label{15}
   X \cdot \psi  = \texttt{i} J(X)\cdot
 \psi,\; \hbox{where} \; X \in TM .
\end{equation}
\begin{lem} If  \ \ $( X +  \texttt{i} Y )  \cdot \psi = 0$ then $ g(X,Y)
=0$ and $g(X,X) = g(Y,Y)$.
\end{lem}
\textbf{Proof}. See the proof in (\cite{Frid}, p. 65) for the
Riemannian case who is valid for the pseudo-Riemannian  case.\\
Lemma 1 implies that $J$ defines an orthogonal  almost complex
structure on $M$. Moreover, from (7) and (15) we obtain $J$ is
parallel, since $\psi$ is parallel. In consequence, $(M,g)$ is a
k$\ddot{a}$hler  manifold.\\
Now if $(M,g)$ is a k$\ddot{a}$hler  manifold, then there exists a
canonical $Spin^c$ structure of $(M,g)$. And
from Remark 1, the following conditions are equivalent \\
(a) $(M,g)$ is not $Spin$, \\
(b)  $H_A = \mathbb{S}^1 $, \\
(c)  $(M,g)$  is not Ricci-flat, \\
(d) $H = U(p',q')$.\\
Then  the equivalence between (i) and (ii) are proved. And from
Theorem 2, we have the  equivalence between (ii) and (iii). To
finish the proof of Theorem 3, it remains to show for $H =
U(p',q')$ that $N=2$. For this, we remark that
 $U(p',q')=SU(p',q') \times U_{\mathbb{S}^1}$, where
$$
U_{\mathbb{S}^1} = \{ \;\; \left(%
\begin{tabular}{llll}
$\lambda $ & 0 & $\cdots $ & 0 \\
0 & 1 & $\ddots $ & $\vdots $ \\
$\vdots $ & $\ddots $ & $\ddots $ & 0 \\
0 & $\cdots $ & 0 & 1
\end{tabular}
\right); \;\;  \lambda \in \mathbb{S}^1 \;\;
 \}
$$
$u(p',q') = su(p',q') \oplus u_{\mathbb{S}^1}$, where
$u_{\mathbb{S}^1} \simeq i\mathbb{R} $ \ \ is the Lie algebra of
$U_{\mathbb{S}^1}$. If we consider the imbedding
$$
\begin{array}{ccc}
  i : u(p',q') & \hookrightarrow & so(2p',2q') \\
  A+ \texttt{i}B =((a_{kl})_{1\leq k ,l \leq n} + \texttt{i}(b_{kl})_{1\leq k ,l \leq n}) &
   \rightarrow & \left(%
   \begin{array}{c} \left(%
\begin{array}{cc}
  a_{kl} & b_{kl} \\
 -b_{kl} & a_{kl} \\
\end{array}%
\right)
\end{array}%
\right)_{1\leq k ,l \leq n} \\
\end{array}
$$
$u_{\mathbb{S}^1}$ is generated by $D_1 = E_{12}$. From \cite{BK},
$u^+ := u(1,...,1)$ and $u^- := u(-1,...,-1)$ generate the space $
V_{su(p',q')}= \{ v \in \Delta_{p,q} \;
 ; \;\; su(p',q') \cdot v=0 \}
 $.  Moreover, by (1)
$$
D_1 \cdot u^+ = \texttt{i} u^+ \; \; \hbox{and}\; \; D_1 \cdot u^-
= -\texttt{i} u^- .
$$
Then $u^+$ and $u^-$ generate the space
$$
 V_{u(p',q')}= \{ v \in \Delta_{p,q} \;
 ; \;\; u(p',q') \cdot v = \texttt{i} \mathbb{R}\; v \}.
 $$
And  the proof of Theorem 3 is finished.

\end{document}